\documentclass[12pt]{amsart}
\usepackage{amscd}

\newtheorem{theorem}{Theorem}
\newtheorem{corollary}{Corollary}

\theoremstyle{definition}

\theoremstyle{remark}
\newtheorem*{remark}{Remark}

\usepackage{amsfonts}
\newcommand{\field}[1]{\mathbb{#1}}
\newcommand{\Q}{\field{Q}}
\newcommand{\R}{\field{R}}
\newcommand{\Z}{\field{Z}}
\newcommand{\C}{\field{C}}
\newcommand{\A}{\field{A}}

\newcommand{\bs}{\backslash}
\newcommand{\ra}{\rightarrow}

\begin{document}

\title[Virtual Betti Numbers of Compact Locally Symmetric Spaces]
{Virtual Betti Numbers  of Compact Locally Symmetric Spaces}

\author{T. N. Venkataramana}

\address{School of Mathematics, Tata Institute of Fundamental 
Research, Homi Bhabha Road, Bombay - 400 005, INDIA.}

\email{venky@math.tifr.res.in}

\subjclass{Primary 11F75; Secondary 22E40, 22E41\\
T. N. Venkataramana,
School of Mathematics, Tata Institute of Fundamental 
Research, Homi Bhabha Road, Bombay - 400 005, INDIA.
venky@math.tifr.res.in}

\date{}

\begin{abstract}
We show that the virtual Betti number of a compact locally symmetric space 
with arithmetic fundamental group is either $0$ or else is infinite.  
\end{abstract}

\maketitle

\section{Introduction}
Let $G$ be a connected non-compact linear Lie group with finite centre, 
such that 
$G$ is simple modulo its centre. Let $\Gamma $ be a torsion free cocompact 
arithmetic (not necessarily congruence ) subgroup in  $G$ and let 
$i\geq 0$ be an integer. Consider the 
direct limit cohomology group 
\[ {\mathcal H}^i= lim H^i(\Delta, \C)\]
where the direct limit is  over all finite index subgroups $\Delta$ in
$\Gamma $; we emphasize that $\Gamma  $ is only assumed to be an arithmetic 
subgroup  of $G$ and is not assumed  to be  a congruence
subgroup of $G$. The dimension of the direct limit ${\mathcal H}^i$ as
a $\C$-vector space is called the {\bf virtual} $i$-th
Betti number of $\Gamma$.
     
\begin{theorem} If the direct limit ${\mathcal H}^i$ is finite 
dimensional, then ${\mathcal  H}^i=H^i(G_u/K,\C)$ where $G_u/K$ is the
compact dual of the symmetric space $G/K$ of $G$.
\end{theorem}

As a special case we recover the following result of Cooper, Long and Reid 
(see \cite{CLR}). 
\begin{corollary} If $M$ is a compact arithmetic hyperbolic $3$-manifold 
with non-vanishing first Betti number, then $M$ has infinite virtual 
first Betti number.   
\end{corollary}
\begin{proof} Take $G=SL_2(\C)$ in Theorem 1, and observe that the 
compact dual $G_u/K=S^3$ has vanishing first cohomology. 
\end{proof}

The present note was motivated by the recent preprint \cite{CLR} of 
Cooper, Long and Reid, where they prove Corollary 1, by using crucially,  
the fact that $M$ is a hyperbolic $3$-manifold. We show that this is true 
in greater generality. The point of Theorem 1 is that the group 
$\Gamma  $ is not assumed to be a congruence subgroup; if $\Gamma $ is a 
congruence subgroup, this is a result of A.Borel (see \cite{B}). 

\section{Proof of Theorem 1}

Let $K\subset G$ be a maximal compact subgroup; write $\mathfrak k$ 
and $\mathfrak g$ for the complexified Lie algebras of $K$ and $G$. 
We have the Cartan decomposition 
${\mathfrak g}={\mathfrak k} \oplus {\mathfrak p}$. Note that $\Gamma $ 
(and hence the finite index subgroup $\Delta$) is torsion-free and 
cocompact in $G$. We then get by the Matsushima-Kuga formula (see \cite{BoW}), 
\[ H^i(\Delta , \C)=Hom _K(\wedge ^i \mathfrak p, 
{\mathcal C}^{\infty}(\Delta \bs G)(0).\]  
In this formula, ${\mathcal C}^{\infty}(\Delta \bs G)(0)$ denotes the 
space of complex valued smooth functions on the manifold $\Delta \bs G$ 
which are annihilated by the Casimir of $\mathfrak g$ (the latter space 
in the Matsushima-Kuga formula may be identified with the space of 
{\bf harmonic} differential forms of degree $i$ on $\Delta \bs G/K$ 
with respect to the $G$-invariant metric on the symmetric space $G/K$). \\ 

Taking direct limits in 
the Matsushima -Kuga formula yields the equality
\[{\mathcal H}^i=lim H^i(\Delta, \C)=Hom _K(\wedge ^i {\mathfrak p}, 
\bigcup _{\Delta \subset 
\Gamma }{\mathcal C}^{\infty}(\Delta \bs G)(0)).\] 
Here, $\Delta $ runs through finite index subgroups of $\Gamma $. 
Consider the space
\[{\mathcal F}=\bigcup _{\Delta \subset \Gamma }{\mathcal \C}^{\infty}(\Delta \bs G)(0).\]
On the space ${\mathcal F}$, $G$ 
acts on the right (since the Casimir commutes with the $G$-action). \\

Now, $\Gamma $ is an arithmetic subgroup of $G$. That is, there is a 
semi-simple (simply connected) algebraic group 
${\bf G}$ defined over $\Q$ and a smooth surjective  homomorphism 
$\pi : {\bf G}(\R)\ra G$ with compact kernel such that $\pi ({\bf G}(\Z))$ 
is commensurable to $\Gamma$.  We define $G(\Q)$ simply to mean the 
image group 
$\pi ({\bf G}(\Q))$. It follows from weak approximation (\cite{PR}) 
that $G(\Q)$ is dense in $G$. \\ 

Now, there is an action on  $\mathcal F$ by $G(\Q)$ on the left (which
therefore commutes  with the  right $G$ action),  as follows.  Given a
function $\phi  \in {\mathcal F}$  and given an element  $g\in G(\Q)$,
the function $\phi $ is left $\Delta $-invariant for some finite index
subgroup  $\Delta  $ in  $\Gamma$.   Consider  the function  $g(\phi)=
x\mapsto  \phi  (g^{-1}x)$.  This  function  is  left-invariant  under
$g\Delta g^{-1}$  and hence under $\Gamma \cap  g\Delta g^{-1}$; since
$g\in G(\Q)$,  it follows that  $g$ commensurates $\Gamma $  and hence
that the subgroup  $\Gamma \cap g\Delta g^{-1}$ is  of finite index in
$\Gamma $.  Therefore, $g(\phi)$ lies in ${\mathcal  F}$. This defines
an action of  $G(\Q)$ on the direct limit  ${\mathcal H}^i$. Note that
under this  action, the  action of $\Delta  $ on the  cohomology group
$H^i(\Delta,\C)$ is trivial. \\

Suppose that ${\mathcal H}^i$  is finite dimensional. Since ${\mathcal
H}^i$  is a  direct  limit  of finite  dimensional  vector spaces,  it
follows that it coincides with  one of them.  Therefore there exists a
finte index subgroup $\Delta $ of $\Gamma $ such that
\[{\mathcal H}^i=H^i(\Delta ,\C).\] 
The last sentence  of the foregoing paragraph says  that while $G(\Q)$
acts on $H^i(\Delta,\C)$, the action by $\Delta$ is trivial. Hence the
action by the normal subgroup $N$ generated by $\Delta $ in $G(\Q)$ is
also trivial.  The  density of $G(\Q)$ in $G$ is  easily seen to imply
the  density of  the normal  subgroup $N$  in $G$.  Thus the  image of
$\wedge ^i{\mathfrak p}$ under any element of ${\mathcal H}^i$ (viewed
via the  Matsushima-Kuga formula as a  (K-equivariant) homomorphism of
$\wedge ^i {\mathfrak p}$ into ${\mathcal F}$), goes into $G$ invariant
functions in  ${\mathcal C}^{\infty}(\Delta\bs G)$,  i,e, the constant
functions. But  $Hom _K(\wedge ^i  {\mathfrak p},\C)$ is the  space of
harmonic  differential  forms on  the  compact  dual  $G_u/K$, and  is
therefore isomorphic to $H^i(G_u/K,\C)$.

This proves Theorem 1. 

\begin{remark}
If $\Gamma $ and all the subgroups $\Delta$ are {\bf congruence} subgroups, 
then one sees at once from strong approximation, that the above 
$G(\Q)$ action on the direct limit translates into the action of 
the ``Hecke Operators'' $G(\A_f)$ ($\A_f$ are the ring of finite adeles) 
and amounts to the proof of Borel in \cite{B}. 
In this sense, the proof of Theorem 1 is an extension of 
Borel's proof to the non-congruence case.   
\end{remark}

\noindent{\bf  Acknowledgement:}   The author gratefully acknowledges 
the hospitality of the Forschungsinstitut f\"ur Mathematik, ETH, 
Zurich in November 2006, where this note was written.


\begin{thebibliography}{JPSH}




\bibitem[CLR]{CLR} D.Cooper, D.D. Long and A.W.Reid, 
On the virtual Betti numbers of 
arithmetic hyperbolic 3-manifolds, preprint October 2,  2006.


\bibitem[B]{B} A.Borel, Cohomologie de sous-groupes discrets et 
representations 
de groupes semi-simples, Asterisque, {\bf 32-33} (1976), 73-112.  


\bibitem[BoW]{BoW} A. Borel and  N. R. Wallach, Continuous cohomology,
discrete subgroups and representations  of reductive groups, Annals of
Math. Studies, {\bf 94}, Princeton Univ. Press, Princeton, 1980.

\bibitem[PR]{PR} V.P.Platonov and A.Rapinchuk, Algebraic Groups and Number Theory, 1990. 




\end{thebibliography}
\end{document}